\documentstyle[12pt, epsfig, epsf, verbatim, amssymb,latexsym]{amsart} 
\newtheorem{thm}{Theorem}[section]

\newtheorem{proposition}[thm]{Proposition}

\newtheorem{dummy}{Theorem}

\newtheorem{dummyc}{Corollary}
\newtheorem{corollary}[thm]{Corollary}
\newtheorem{lem}[thm]{Lemma}
\theoremstyle{definition}
\newtheorem{definition}[thm]{Definition}

\theoremstyle{remark}
\newtheorem{example}[thm]{$\mathrm{Example}$}
\newtheorem{remark}[thm]{$\mathrm{Remark}$}

\def\u{{\rm{u}}}
\def\d{{\rm{d}}}
\def\bu{{\rm{u}}_{\rm{b}}}
\def\bd{{\rm{d}}_{\rm{b}}}
\def\la{\rho}
\def\Z{\text{\boldmath $Z$}}

\def\Arf{{\rm{Arf}}}
\def\sign{{\rm{sign}}}
\def\Ms{\sigma}
\def\om{\omega}

\numberwithin{equation}{section}
\begin{document}

\title[Unoriented band surgery on knots and links]
{Unoriented band surgery on knots and links}

\date{\today}

\author{Tetsuya Abe and Taizo Kanenobu}
\address{Research Institute for Mathematical Sciences, 
Kyoto University, Kyoto 606-8502 Japan} 
\email{tetsuya@@kurims.kyoto-u.ac.jp}

\address{Department of Mathematics, Osaka City University,
Sugimoto, Sumiyoshi-ku,
Osaka 558-8585, Japan} 
\email{kanenobu@@sci.osaka-cu.ac.jp}



\maketitle

\begin{abstract}
We consider a relation between two kinds of unknotting numbers defined by using a band surgery on unoriented knots; the band-unknotting number and $H(2)$-unknotting number, which we may characterize in terms of the first Betti number of surfaces in $S^3$ spanning the knot and the trivial knot.
We also give several examples for these numbers.
\end{abstract}

\section{Introduction}
\begin{figure}[htbp]
\centerline{\includegraphics[scale=.5]{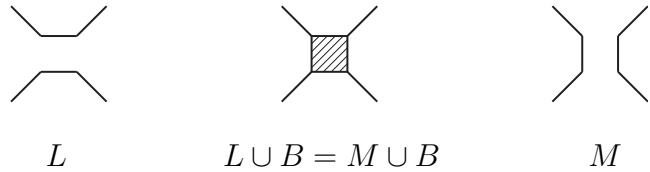} }
\par \vspace{5pt}
\centerline{$L$ \hspace{18mm} $L\cup B=M \cup B$  \hspace{17mm} $M$}
\caption{$M$ is obtained form $L$ by a band surgery along the band $B$, and vice virsa.} 
\label{bandsurgery}
\end{figure}

A band surgery is an operation which deforms a link into another link.
Let $L$ be a link and  $b : I \times I \to S^3$ 
an embedding 
such that $L \cap b(I \times I) = b(I \times \partial I)$, 
where $I$ is the unit interval $[0, 1]$. 
Then  we may obtain a new link
$M= \left (  L \setminus  b(I \times \partial I) \right ) \cup b(\partial I \times I)$,
which is 
called a link obtained from $L$ 
by the \emph{band surgery} along the band $B$, where  $B=b(I \times I)$; 
see Fig.~\ref{bandsurgery}.

A band surgery appears in various aspects in knot theory.
For example, it is an essential tool in the study of surfaces embedded in $4$-space; 
it is  motivated by the study of DNA site-specific recombinations; 
cf.~\cite{MR1068451, ShimokawaKoya:2009-04-20}.
In this paper, we study an unoriented band surgery, that is, we consider a band surgery for unoriented knots and links.
Using an unoriented band surgery,  we can define two numerical invariants for a knot.
One is the \emph{band-unknotting number} of a knot $K$,   $\bu(K)$, which is the minimal number of band surgeries to deform $K$ into the trivial knot.
The other is the  $H(2)$-\emph{unknotting number} of $K$,  $\u_{2}(K)$, 
introduced by Hoste, Nakanishi, and Taniyama \cite{HNT90},
which is the minimal number of component-preserving band surgeries to deform  $K$ into the trivial knot, and so the 
unknotting sequence which realizes the $H(2)$-unknotting number consists of only knots, whereas the unknotting sequence which realizes the band-unknotting number may contain a link with more than one component.
Then by definition $\bu(K) \le  \u_2(K)$.
Moreover, we have:
\renewcommand{\thedummyc}{\ref{cor:yasuhara2}}
\begin{dummyc}
For a knot $K$, 
\begin{equation*}
\bu(K) =\u_{2}(K) -1  \,{\text{ or }}\,  
\u_2(K).
\end{equation*}
Furthermore, if $\bu(K)$ is odd, 
then
$\bu(K)=  \u_{2}(K)$; equivalently,
if $\u_2(K)$ is either one or even, 
then
$\bu(K)=  \u_{2}(K)$.
\end{dummyc}
Actually, we have  $\bu(8_{18})=2$ and $\u_{2}(8_{18})=3$ (Example~\ref{ex;8_18}(i)),
which may be generalized as follows:
\renewcommand{\thedummy}{\ref{them;main-theorem}}
\begin{dummy} 
There exist infinitely many knots $K$ such that $\bu(K)=2$ {and} $\u_{2}(K)=3$.
\end{dummy}

On the other hand, Taniyama and Yasuhara \cite[Theorem~5.1]{TaniyamaYasuhara94} characterized the band-unknotting number  $\bu(K)$ as $\min_F \{ \beta_1(F)\} -1$, where the minimum is taken over all connected (orientable or nonorientable)  surfaces $F$ in  $S^3$ spanning the knot $K$ and the trivial knot, where $ \beta_1(F)$ is the first Betti number of $F$. 
Similarly, we may characterize the $H(2)$-unknotting number ${\u}_{2} (K)$ as $\min_{F_N} \{ \beta_1(F_N)\} -1$, where the minimum is taken over all connected  \emph{nonorientable}  surfaces $F_N$ in  $S^3$ spanning the knot $K$ and  the trivial knot (Theorem~\ref{thm:yasuhara4}).
In the following, we will consider 
the band-Gordian distance and $H(2)$-Gordian distance,
which generalize the band-unknotting number and  $H(2)$-unknotting number, respectively.

This paper is organized as follows:
In Sec.~\ref{sect;band Gordian distance}, 
we give definitions of band-Gordian distance and band-unknotting number, and
those of $H(2)$-Gordian distance and  $H(2)$--unknotting number.
In Sec.~\ref{sect;band-surgeries and surfaces in S^3}, we prove Theorem~\ref{thm:yasuhara2}, which generalizes Corollary~\ref{cor:yasuhara2} to the distance.
In Sec.~\ref{sect;characterization}, we give  characterizations of the band-Gordian distance and $H(2)$-Gordian distance in terms of the first Betti number of surfaces in $S^3$ spanning the knots  (Eq.~\eqref{eq;TYTheorem5.1} and Theorem~\ref{thm:yasuhara4}).
In Sec.~\ref{section:invariants}, 
we review some invariants; 
the homology groups of the cyclic branched covering space of a link, the Jones 
polynomial, the Arf invariant, the signature, and the Q polynomial.
We use them to evaluate the band-Gordian distance. 
In Sec.~\ref{sect;Relations to the band- and usual Gordian distances},
we  give estimations for the band-Gordian distance and 
the band-unknotting number in terms of the usual Gordian distance and unknotting number, respecitevely (Theorem~\ref{thm;u1knot}).
Using Theorem~\ref{thm;u1knot} and invariants given in Sec.~\ref{section:invariants} 
we determine the band-Gordian distances and band-unknotting numbers for examples with $H(2)$-Gordian distance or $H(2)$-unknotting number 3 given in \cite{Kanenobu11}.  
In Sec.~\ref{sect;example}, we prove Theorem \ref{them;main-theorem}.
In Sec.~\ref{section:9cr-knots},
we give a complete table of band-unknotting numbers and $H(2)$-unknotting numbers for knots with up to $9$ crossings.
In Sec.~\ref{sect:Special values of the Jones polynomial and a slice knot}, we give a relation for a special value of the Jones polynomial of a  slice knot.
\section{Band-Gordian distance and $H(2)$-Gordian distance} 
\label{sect;band Gordian distance}

\begin{figure}[htbp]
\centerline{ \includegraphics[scale=.5]{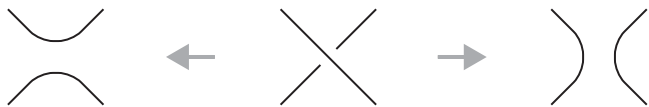} }
\caption{Smoothing a crossing.} 
\label{smoothing}
\end{figure}

In this section,
we introduce the band-Gordian distance and $H(2)$-Gordian distance.
It is easy to see that any connected unoriented link diagram can be deformed into the trivial knot diagram without any crossing by smoothing appropriately at every crossing as shown in Fig.~\ref{smoothing}; cf.~\cite[Part~1, Sec.~4]{Kauffman93}.
Since a smoothing is realized by an unoriented band surgery, any link can be deformed into the trivial knot by a sequence of unoriented band surgeries.
Thus any two links are related by a sequence of unoriented band surgeries.

\begin{definition}
Let $L$ and $M$ be unoriented links.
The \emph{band-Gordian distance} from  $L$  to $M$,
 denoted by $\d_{b}(L, M)$,  is defined to be 
the minimal number of unoriented band surgeries needed to deform $L$ into $M$.
In particular, the  \emph{band-unknotting number} of a knot $K$,
 $\bu(K)$,  is defined to be the band-Gordian distance from  $K$ to the trivial knot $U$;
 $\bu(K)=\bd(K, U)$.
\end{definition}

Hoste, Nakanishi, and Taniyama \cite{HNT90} introduced an $H(n)$-move, which is a deformation of a link diagram.  In particular, an \emph{$H(2)$-move} is a band surgery which requires to keep the number of the components. 
An $H(2)$-move is an unknotting operation, that is, any knot $K$ can be deformed into the trivial knot by a sequence of $H(2)$-moves \cite[Theorem~1]{HNT90}.
Let $J$ and $K$ be knots.
The  \emph{$H(2)$-Gordian distance} from $J$ to $K$, denoted by $\d_2(J,K)$, is defined to be the minimal number of $H(2)$-moves 
needed to transform $J$  into $K$, and the  \emph{$H(2)$-unknotting number} of a knot $K$,  denoted by $\u_2(K)$,  is defined to be the $H(2)$-Gordian distance from  $K$ to the trivial knot $U$;
 $\u_2(K)=\d_2(K, U)$.

\section{Band surgeries and surfaces in $S^3$} 
\label{sect;band-surgeries and surfaces in S^3}
For two knots $J$ and $K$, by definition we have 
$\bd(J,K)\le \d_2(J,K)$.  
Moreover we have:
\begin{thm} \label{thm:yasuhara2}
For any knots $J$ and $K$, we have 
\begin{equation}
\bd(J,K) =\d_{2}(J,K) -1 \,{\text{ or }}\,  \d_2(J,K).
\label{eq;thm:yasuhara2}\end{equation}
Furthermore, if $\bd(J,K)$ is odd, 
then
$\bd(J,K)=  \d_{2}(J,K)$; equivalently,
if $\d_2(J,K)$ is either one or even, 
then
$\bd(J,K)=  \d_{2}(J,K)$.
\end{thm}
In particular, we have:
\begin{corollary} \label{cor:yasuhara2}
For a knot $K$, 
\begin{equation}
\bu(K) =\u_{2}(K) -1 
\,{\text{ or }}\,  \u_2(K).
\end{equation}
Furthermore, if $\bu(K)$ is odd, 
then
$\bu(K)=  \u_{2}(K)$; equivalently,
if $\u_2(K)$ is either one or even, 
then
$\bu(K)=  \u_{2}(K)$.
%
\end{corollary}

We prove Theorem~\ref{thm:yasuhara2} using Lemma~\ref{lemma:yasuhara3} below, which allows us to understand  the band-Gordian distance $\bd(J,K)$ in terms of the first Betti number of a surface in $S^3$ bounding 
a  2-component link whose components are isotopic to $J$ and $K$.
Although  Lemma~\ref{lemma:yasuhara3} is essentially given by Taniyama and Yasuhara  \cite[Theorem~5.1]{TaniyamaYasuhara94}, for the sake of completeness we provide a proof.


We may perform a band surgery along a several number of disjoint bands simultaneously.
Let $L$ be a link and $b_i: I \times I \to S^3$ ($i=1, \dots, n$) are disjoint embeddings, that is,
$b_{i}(I \times I) \cap b_{j}(I \times I) = \emptyset$ for $i \neq j$ such that
$L \cap b_{i}(I \times I) = b_{i}(I \times \partial I)$.
Then we obtain another link $M$;
\begin{equation}
M=\left( L \setminus \bigcup_{i=1}^n  b_i(I \times \partial I) \right) \cup \bigcup_{i=1}^n b_i (\partial I \times I).
\end{equation}
Then $M$ is called the link obtained from $L$ by the (multiple) band surgery along the bands $B_1, \dots , B_n$, where $B_i=b_i(I\times I)$, which we denote by $h(L; B_1, \dots, B_n)$.

In order to give a nomal form for an embedded surface in $S^4$ the following lemma is essentially proved in {\cite[Lemma 1.14]{MR672939}} and {\cite[Lemma 2.4]{Kamada89}}.
\begin{lem}\label{lem:Kawauchi-Shibuya-Suzuki82}
Suppose that $L$ and $M$ are links, which are related by a sequence of $n$ band surgeries.
Then there exist mutually disjoint bands $B_1, \dots, B_n$ such that
$M=h(L; B_1, \dots, B_n)$.
\end{lem}

\begin{pf}
We use induction on $n$.
The case $n=1$ is trivial.
Let $L_{i}$, $i=0,1, \dots, n$, be links 
such that
$L_{0}=L$,  $L_n=M$,
and $L_{i}$ is obtained from $L_{i-1}$ by a band surgery along $b_i$;
$L_{i}=h(L_{i-1}, B_{i})$, $1\le i  \le n$,  where $B_{j}=b_{j}(I \times I)$.
Suppose that the lemma holds for $n-1$. 
Then there exist mutually disjoint  bands 
$\tilde{B}_{1}, \dots , \tilde{B}_{n-1}$ 
such that 
$L_{n-1}=h(L; \tilde{B}_{1}, \dots , \tilde{B}_{n-1})$.
Let $\alpha$ be an attaching arc of $B_{n}$;
$\alpha = b_{n}(I\times \{t\})$,  $t=0$, $1$. 
If 
$\alpha \cap \tilde{B}_{j} \neq \emptyset$, $j<n$,
we slide $B_n$ along $L_{n-1}$ as shown in Fig.~\ref{band-slide}.
Therefore we may assume that
$\alpha \cap (\tilde{B}_{1} \cup \cdots \cup \tilde{B}_{n-1}) = \emptyset$.

\begin{figure}[htbp]
\centerline{\includegraphics[scale=.8]{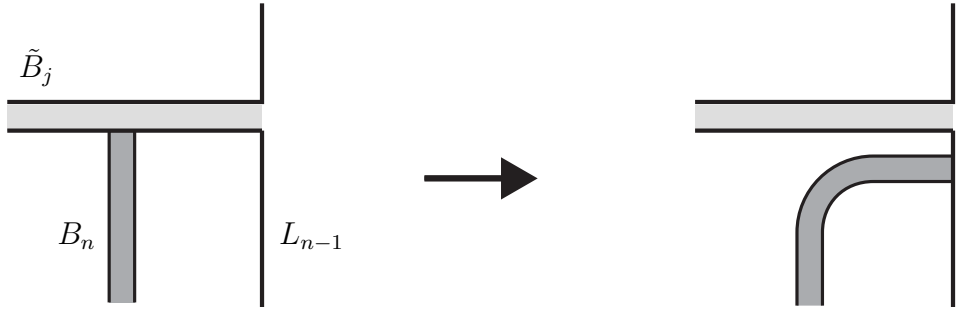} }
\vspace*{-37mm}\hspace*{-120mm}$\tilde B_j$\\
\vspace*{45pt}\hspace*{-75mm}$B_n$ \hspace*{23mm}$L_{n-1}$
\vspace*{15pt}
\caption{Sliding the band $B_n$ along $L_{n-1}$.} 
\label{band-slide}
\end{figure}

Since we may assume that
$B_{n}$ and $\tilde{B}_{1} \cup \cdots \cup \tilde{B}_{n-1}$ 
intersect transversely,
the intersection $B_{n} \cap (\tilde{B}_{1} \cup \cdots \cup \tilde{B}_{n-1})$
is a 1-manifold,  that is, simple loops and arcs.
Now we choose a proper simple arc $\gamma$ in  $B_{n}$
connecting the two attaching arcs of $B_{n}$ 
as in Fig.~\ref{band-tube}(a).
Replace $B_{n}$ with a sufficiently small regular neighborhood of $\gamma$
in $B_{n}$.
Then we may assume that the intersection
$B_{n} \cap (\tilde{B}_{1} \cup \dots \cup \tilde{B}_{n-1})$ 
consists of proper simple arcs in $B_{n}$ as in  Fig.~\ref{band-tube}(b). 
In particular, 
$B_{n}  \cap (\partial \tilde{B}_{1} \cup \dots \cup \partial  \tilde{B}_{n-1}) = \emptyset$.

\begin{figure}[htbp]
{\centerline{\includegraphics[scale=.7]{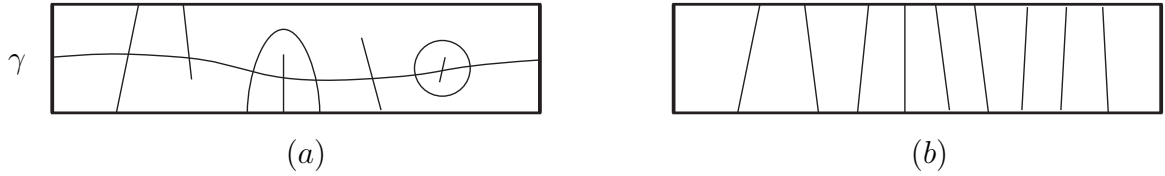} }}
\par\vspace*{-38pt}\hspace*{-450pt}$\gamma$\par %
\vspace*{18pt}
\par\centerline{ $(a)$ \hspace{75mm} $(b)$}
\caption{Repalacing the band $B_n$.} 
\label{band-tube}
\end{figure}

Finally, we deform the band $B_{n}$ so that 
$B_{n}  \cap (\tilde{B}_{1} \cup \dots \cup   \tilde{B}_{n-1}) = \emptyset$ as in Fig.~\ref{band-slide2}.
Let $\tilde{B}_{n}$  be the resulting band.
Then the bands $\tilde{B}_{1}, \dots ,\tilde{B}_{n}$ are 
mutually disjoint and
$M=h(L; \tilde{B}_{1}, \dots , \tilde{B}_{n})$.
Therefore the lemma holds for $n$, completing the proof.
\end{pf}

\begin{figure}[htbp]
{\centerline{\includegraphics[scale=.7]{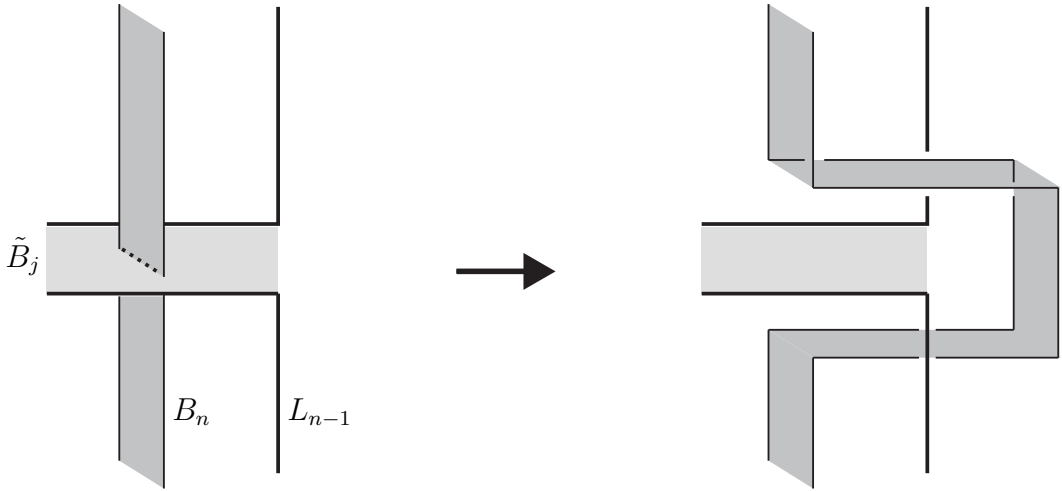} }}
\vspace*{-105pt}\hspace*{-143mm}$\tilde B_j$\\
\vspace*{40pt}\hspace*{-78mm}$B_n$ \hspace*{25pt}$L_{n-1}$
\vspace*{20pt}
\caption{Deforming the band $B_n$.} 
\label{band-slide2}
\end{figure}

\begin{lem} \label{lemma:yasuhara3}
If two knots $J$ and $K$ are related by a sequence of $n$ band surgeries,
then there exists a compact connected surface $F$ in $S^3$ 
such that:\\ 
{\rm(i)}  $F$ spans the two knots $J$ and $K$, that is, 
$\partial F$ is a $2$-component link $\tilde{J} \cup \tilde{K}$ with $\tilde J$ and $\tilde{K}$  isotopic to $J$ and $K$, respectively; and \\
{\rm(ii)} $\beta _{1} (F)=n+1$,\\ 
where $\beta _{1} (F)$ is the first Betti number of $F$.
\end{lem}

\begin{pf}
By Lemma~\ref{lem:Kawauchi-Shibuya-Suzuki82}
there exist mutually disjoint  bands $B_{1}, \dots , B_{n}$ such that 
$K=h(J; {B}_{1}, \dots , {B}_{n})$.
We take an annulus $A$ embedded in $S^3$ so that 
one of the boundary components is $J$ and 
$A\cap B_i=J\cap B_i$ is the attaching arcs of $B_i$ for each $i$.
Let $F$ be the surface  $A \cup \ (\bigcup_{i=1}^{n} {B}_{i})$.
Then $F$ satisfies the conditions (i) and (ii).
\end{pf}


\begin{pf*}{Proof of Theorem~$\ref{thm:yasuhara2}$}
By definition, 
we have $\bd(J,K) \le \d_{2}(J,K)$.  
Let $n=\bd(J,K)$.
By Lemma \ref{lemma:yasuhara3},
there exists a compact connected surface $F$ such that
$\partial F= \tilde{J} \cup \tilde{K}$ and $\beta_1(F)=n+1$, 
where $\tilde{J}$ and $\tilde{K}$ are ambient isotopic to $J$ and $K$,  respectively.
There are two cases.

(i) $F$ is nonorientable. 
We can deform $F$ into the surface as shown in Fig.~\ref{nonorientable}(a),
which is constructed by attaching $n$ half-twisted bands and $1$ bands to a disc; 
cf.~\cite{Clark78}.  In fact,  $F$ is the connected sum of $n$ projective planes minus two open discs, which is homeomorphic to the surface as shown in  Fig.~\ref{nonorientable}(b).
Then adding  $n$ bands to $\tilde J$ as shown in  Fig.~\ref{nonorientable}(c), 
we obtain a parallel link which bounds an annulus whose center line is isotopic to $K$. 
This implies that the knot $J$ is transformed into $K$ by performing the $H(2)$-move $n$ times. 
Thus $\d_2(J,K)\le n$, and so we have $\bd(J,K) =\d_{2}(J,K)$.

\begin{figure}[htbp]
{\centerline{\includegraphics[scale=.5]{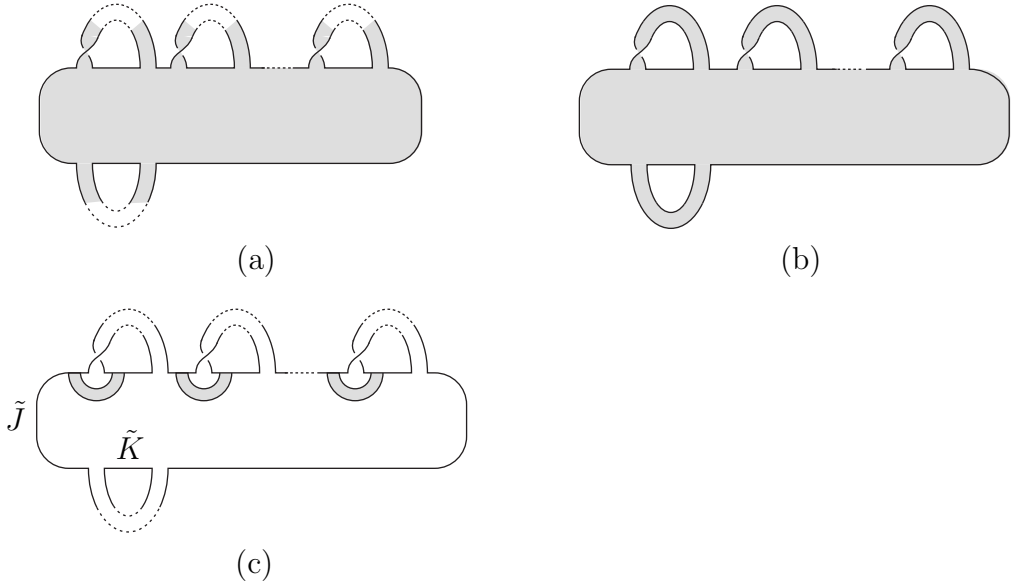} }}
\par\vspace*{-42mm}(a)\hspace*{67mm}(b)\\
\vspace*{15mm}\hspace*{-138mm} $\tilde{J}$\\
\vspace*{-7pt}\hspace*{-106mm}$\tilde{K}$
\par\vspace{25pt}(c)\hspace*{73mm}
\caption{The case $F$ is nonorientable.} 
\label{nonorientable}
\end{figure}

(ii) $F$ is orientable.  
Then $F$ is a Seifert surface for the link $\tilde{J} \cup \tilde{K}$ and is represented in the disc-band form
as shown in Fig.~\ref{orientable}(a), which is constructed by attaching $(n+1)$ bands to a disc,
and so $n$ is an even number and the genus of $F$ is $n/2$; cf.~\cite{MR0052104}, 
Attaching a half-twisted band to $F$, we obtain a nonorientable surface $F'$ 
as shown in Fig.~\ref{orientable}(b).  
Then $\beta _{1} (F')=n+2$ and the boundary of $F'$ is isotopic to the link $\tilde{J} \cup \tilde{K}$.
Therefore, from Case (i) we have $\d_2(J,K)\le n+1$, giving Eq.~\eqref{eq;thm:yasuhara2}.
Noting that $n$ is even, the proof is complete.
\end{pf*}

\begin{figure}[htbp]
{\centerline{\includegraphics[scale=.5]{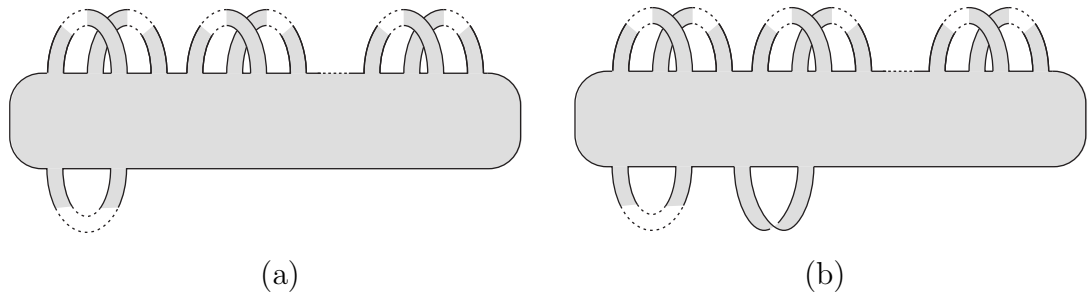} }}
\par(a)\hspace*{67mm}(b)
\caption{The case $F$ is orientable.} 
\label{orientable}
\end{figure}

\section{Characterizations of the band-Gordian distance and $H(2)$-Gordian
distance}
\label{sect;characterization}
For two knots $J$ and $K$, Taniyama and Yasuhara \cite[Sec.~5]{TaniyamaYasuhara94} defined two indices.
\begin{itemize}
\item $\tilde{\d}_C (J,K)=\min_F \{ \beta_1(F)\} -1$, where the minimum is taken over all connected (orientable or nonorientable)  surfaces $F$ in  $S^3$ spanning the two knots $J$ and $K$.
\item $\tilde{c}(J,K)$, the minimum number of critical points of a locally flat (orientable or nonorientable) surface in $S^3\times [0,1]$ bounded by $J\times \{0\}$ and $K\times\{1\}$.
\end{itemize}
Then they showed \cite[Theorem~5.1]{TaniyamaYasuhara94}:
\begin{equation}
\bd(J,K)=\tilde{\d}_C (J,K)=\tilde{c}(J,K).
\label{eq;TYTheorem5.1}
\end{equation}
In fact, from the proof of Theorem~\ref{thm:yasuhara2} we see the first equality of Eq.~\eqref{eq;TYTheorem5.1}.

\begin{remark}
In \cite{TaniyamaYasuhara94}, the symbol $\d_2(J,K)$ is used to mean the band-Gordian distance $\bd(J,K)$, but not our $H(2)$-Gordian distance.
\end{remark}

We define two analogous indices:
\begin{itemize}
\item $\tilde{\d}_{C_N} (J,K)=\min_{F_N} \{ \beta_1(F_N)\} -1$, where the minimum is taken over all connected  \emph{nonorientable}  surfaces $F_N$ in  $S^3$ spanning two knots $J$ and $K$.
\item $\tilde{c}_N(J,K)$, the minimum number of critical points of a locally flat \emph{nonorientable} surface in $S^3\times [0,1]$ bounded by $J\times \{0\}$ and $K\times\{1\}$.
\end{itemize}
Then we have:
\begin{thm} \label{thm:yasuhara4}
\begin{equation}
\d_2(J,K)=\tilde{\d}_{C_N} (J,K)=\tilde{c}_N(J,K).
\end{equation}
\end{thm}

\begin{pf} We only prove the first equality.
Let $\d_2(J,K)=n$.
By Lemma~\ref{lemma:yasuhara3} there exists a nonorientable surface $F$ spanning two knots $J$ and $K$ with $\beta_1(F)=n+1$.  In fact, let $K_0$, $K_1, \dots, K_n$ be a sequence of knots such that $K_0 = J$, $K_n = K$, and $K_i$ is obtained from $K_{i-1}$ by  a single $H(2)$-move.  
Since $K_1$ is obtained from $J$ by adding a half-twisted band,  the surface $F$ constructed as in the proof of Lemma~\ref{lemma:yasuhara3} is nonorientable.
Therefore, $ \tilde{\d}_{C_N} (J,K)\le n$.

Conversely, let $n=\tilde{\d}_{C_N}(J,K)$,  and let $F$ be a nonorientable surface spanning  the two knots $J$ and $K$ with $\beta_1(F)=n+1$.  
Then $F$ may be illustrated as in Fig.~\ref{nonorientable}(a).  By adding $n$
 bands as in  Fig.~\ref{nonorientable}(c), $J$ is transformed into $K$, and so we have  $\d_2(J,K) \le n$, completing the proof.
\end{pf}

\section{Some invariants}
\label{section:invariants}

In this section, 
we review the definitions and some properties of the homology groups of the cyclic branched covering space,
 Jones polynomial, Arf invariant,  signature,  and Q polynomial,
 which enable us to estimate the band-Gordian distance and band-unknotting number.

\subsection{The homology group of the cyclic branched covering space} 
\label{subsect;cover}

For an oriented link $L$, let $\Sigma_p(L)$ be the $p$-fold cyclic covering space of $S^3$ branched along $L$, and
 $e_p(L)$  the minimum number of generators of $H_1(\Sigma_p(L); \Z)$.
If either $p=2$ or $L$ is a knot, then the space $\Sigma_p(L)$ is irrelevant to the orientations of $L$.  Therefore, we can prove the following in a similar way to Theorem~4 in \cite{HNT90}.
\begin{lem} \label{lem;mg} 
{\rm{(i)}} If two links $L$ and $M$ are related by a single unoriented  band surgery, then
$|e_2(L)-e_2(M)|\le 1$.

{\rm{(ii)}}  If two knots $J$ and $K$ are related by a single $H(2)$-move, then
$|e_p(J)-e_p(K)|\le p-1$.
\end{lem}
This immediately implies:
\begin{proposition}\label{prop;mg} 
{\rm{(i)}}  For two links $L$ and $M$ we have: 
\begin{gather}
\bd(L,M)\ge |e_2(L)-e_2(M)|; \\
\bu(L)\ge e_2(L).
\label{eq;mu}
\end{gather}

{\rm{(ii)}} For two knots $J$ and $K$ we have: 
\begin{gather}
\d_2(J,K)\ge |e_p(J)-e_p(K)|/(p-1).
\label{eq;mu2}
\end{gather}
\end{proposition}

\begin{remark}
Proposition~\ref{prop;mg}(i) gives  the same statement as Proposition~5.2  with $p=2$ in \cite{HNT90}, which does not seem to hold for $p>2$.
Proposition~\ref{prop;mg}(ii) with $J$ unknot is given in Theorem~2.1 in \cite{Kanenobu-Miyazawa09}. 
\end{remark}

\begin{example} \label{ex;Kn1}
For each positive integer $i$ let $J_i$ be either the right-hand trefoil knot $3_1$ or the left-hand trefoil knot $3_1!$,
and $K_n$ the connected sum  $J_1\# \cdots\# J_n$.
Then since $H_1(\Sigma_2(J_i); \Z)$ is isomorphic to $\Z_3$,
$H_1(\Sigma_2(K_n);\Z)$ is the $n$-fold direct sum of $\Z_3$.
Thus by Proposition~\ref{prop;mg} we have
$\bd(K_m, K_n)\ge |m-n|$ and $\bu(K_n)\ge n$.
Since $\u_2(3_1)=1$, we obtain
\begin{gather}
\bd(K_m, K_n)= |m-n|;   \label{eq;Kn11} \\
\bu(K_n)=n.
\label{eq;Kn12}\end{gather}
 \end{example}

\subsection{The Jones polynomial} 
\label{subsect;Jones}

\begin{figure}[htbp]
\centerline{\includegraphics[scale=.5]{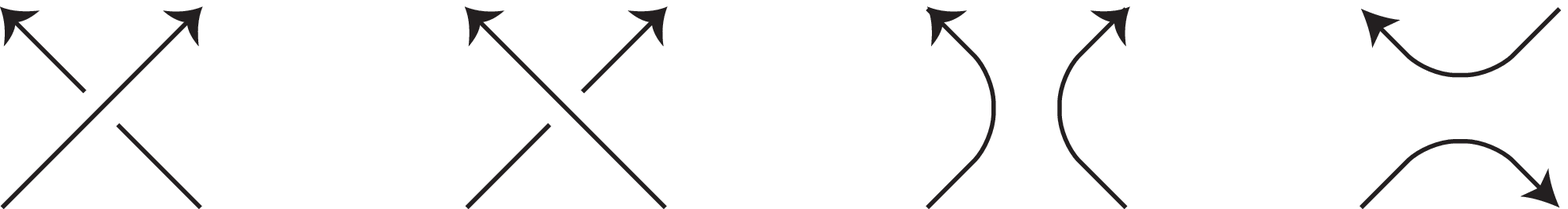}}
\centerline{ $L_+$ \hspace{22mm} $L_-$ \hspace{22mm}  $L_0$ \;}
\caption{A skein triple.} 
\label{skein}
\end{figure}

The \emph{Jones polynomial}  $V(L; t) \in {\Z}[t^{-1/2}, t^{1/2}]$  \cite{Jones87}
of an oriented link $L$ is an isotopy invariant of an oriented link defined by the following formulas:
\begin{gather}
V(U; t) =1; \\
t^{-1} V(L_{+}; t) - t  V(L_{-}; t) 
                = \left( t^{1/2}-t^{-1/2} \right)  V(L_{0}; t) ; 
\label{eq;V_skein}
\end{gather}
where $U$ is the unknot and 
$L_{+}$, $L_{-}$, $L_{0}$  are three links 
that are identical except 
near one point where they are as in  Figure~\ref{skein}.
We call the triple $(L_{+}, L_{-}, L_{0})$ a \emph{skein triple}.  

Let $L$ be a link with $c$ components and
$\delta(L)=\dim H_1(\Sigma_2(L); \Z_3)$.
Then Lickorish and Millett \cite[Theorem~3]{MR860127} 
have shown
\begin{equation}
V(L ; \om) = \pm i^{c-1} (i\sqrt{3})^{\delta(L)},
\label{eq;LickorishMillett86}\end{equation}
where $\om=e^{i\pi /3}$ and $V(L ; \om)$ means the value of $V(L;t)$ at $t^{1/2}=e^{i\pi/6}$.
Notice that if $L'$ is a link obtained from $L$ by changing the orientation of one component of $L$, say $K$, then
\begin{equation}
V(L';t) =  t^{-3\lambda} V(L;t ),
\end{equation}
where $\lambda$ is the  linking number of $K$ and  the remainder of $L$; cf.~\cite[Theorem 11.2.9]{Murasugi96}.  Thus, $V(L ; \om) =  \pm V(L' ;\om)$.
This special value of the Jones polynomial evaluates the band-Gordian distance and band-unknotting number.  
\begin{thm}\label{thm;Jpoly}
If two links $L$ and $M$ are related by a single band surgery, then
\begin{equation}
|V(L ; \om) /V(M ; \om)| \in \{ 1, \sqrt3^{\pm 1}\},
\label{eq;prop;Jpoly1}
\end{equation}
where the orientations of $L$ and $M$ are irrelevant.  
 \end{thm}
This immediately implies the following:
\begin{corollary}\label{coro;Jpoly}  Let $L$ and $M$ be links. 
\par
{\rm(i)} 
It holds that $\bd (L,M) \ge |\delta(L)-\delta(M)|$.  In particular,
$\bd (L) \ge \delta(L)$.

{\rm(ii)} If $|V(L ; \om) /V(M ; \om) |= \sqrt3^n$, 
then $\bd(L, M) \ge |n|$.   In particular,  if $|V(L ; \om)|= \sqrt3^n$,
then $ \bu (L)\ge n$.
\end{corollary}
In order to prove Theorem~\ref{thm;Jpoly} we use the following, which is due to Miyazawa \cite[Proposition~4.2]{Miyazawa11}; cf. \cite[Lemma 2.1]{Kanenobu10}.
\begin{lem} \label{lemma;Miyazawa2010}
Let $(L_{+}, L_{-}, L_{0})$ be a skein triple. Then
\[  V(L_+ ; \om)/ V(L_- ; \om)  \in  \{\pm 1, i \sqrt{3}^{\pm1} \}. \]
\end{lem} 

\begin{pf*}{Proof of Theorem~$\ref{thm;Jpoly}$}
If $L$ and $M$ are related by a single orientable band surgery, then by Theorem 2.2 in \cite{Kanenobu10}
we have
\begin{equation}
V(L ; \om)/ V(M ; \om) \in \{ \pm i,  - \sqrt{3}^{\pm1} \}.
\end{equation}

We consider the case $L$ and $M$ are related by a nonorientable band surgery.
Let $c(L)$ denote the number of components of a link $L$.
For a skein triple $(L_+, L_-, L_0)$, let $L_{\infty}$ be an oriented link
represented by one of the diagrams in Figure~\ref{skein-inf}, where
(i) $c(L_+)<c(L_0)$ and (ii) $c(L_+)>c(L_0)$.
Then we may assume $L_0$ and $L_{\infty}$ are isotopic to $L$ and $M$, respectively.
\begin{figure}[htbp]
\centerline{\includegraphics[scale=.5]{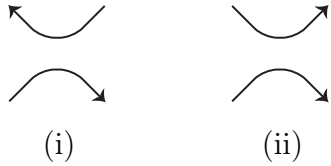}}
\centerline{ (i) \hspace{22mm} (ii)}
\caption{Two possible orientations for $L_{\infty}$.} 
\label{skein-inf}
\end{figure}

Then we have the  `$V_{\infty}$' formuals \cite{Birman88}; cf.~\cite{Lickorish86}:
\begin{align}
 V(L_{+}; t) - t  V(L_{-}; t) +t^{3\lambda}(t-1)  V(L_{\infty}; t)&=0; \label{eq;Vinf1}\\
  V(L_{+}; t) - t  V(L_{-}; t) +t^{3(\nu-\frac12)}(t-1)  V(L_{\infty}; t)&=0, \label{eq;Vinf2}
\end{align}
where Eq.~\eqref{eq;Vinf1} holds for Case (i) and $\lambda$ is the linking number of the right-hand component of $L_0$ in Figure \ref{skein} with the remainder of $L_0$, and
Eq.~\eqref{eq;Vinf2} holds for Case (ii) and $\nu$ is the linking number of the bottom-right to top-left component of $L_+$ in Figure~\ref{skein} with the remainder of $L_+$.
From Eqs.~\eqref{eq;V_skein}, \eqref{eq;Vinf1} and \eqref{eq;Vinf2} we have
\begin{align}
 \dfrac{V(L_0; t)}{ V(L_{\infty}; t)}
 &
 =
 \begin{cases}
  \dfrac{ t^{-{1}/{2} }V(L_+; t) - t^{3/2} V(L_-; t) }{ -t^{-3\lambda} ( V(L_+; t) - t V(L_-; t))} &{\text{for Case (i)}};
  \\
   \dfrac{ t^{-{1}/{2} }V(L_+; t) - t^{3/2} V(L_-; t) }{ -t^{-3\nu+3/2} ( V(L_+; t) - t V(L_-; t))} &{\text{for Case (ii)}}.
\end{cases}
\end{align}
Then letting $x=V(L_+ ; \om)/ V(L_-; \om)$, 
we have
\begin{align}
\dfrac{V(L_0 ; \om)}{ V(L_{\infty} ; \om)} &
=
\begin{cases}
\dfrac{\om^{-1/2}x-\om^{3/2} }{-\om^{-3\lambda}(x-\om )}
=\pm \dfrac{x-\om^{2} }{\om^{1/2}(x-\om)} & {\text{for Case (i)}};
\\
\dfrac{\om^{-1/2}x-\om^{3/2} }{-\om^{-3(\nu-1/2)}(x-\om )}
=\pm i\dfrac{x-\om^{2} }{\om^{1/2}(x-\om)}  &{\text{for Case (ii)}}.
\end{cases}\end{align}
For  $x=1$, $-1$, $-i\sqrt 3$, $-i\sqrt 3^{-1}$, we obtain ${(x-\om^{2}) }/{\om^{1/2}(x-\om)}=\sqrt3$, $\sqrt3^{-1}$, $-i$, $i$, respectively.
Thus by Lemma~\ref{lemma;Miyazawa2010}  the proof is complete.
\end{pf*}

\begin{example} \label{ex;Kn2}
We use the same notation as in Example~\ref{ex;Kn1}. 
Since $V(3_1; \om)=-i\sqrt 3$ and $V(3_1!;\om)=i\sqrt 3$,
we have $|V(K_n; \om)|={\sqrt 3}^n$, 
and so  Corollary~\ref{coro;Jpoly}(ii)  also implies  Eqs.~\eqref{eq;Kn11} and \eqref{eq;Kn12}.
 \end{example}

\subsection{The Arf invariant} 
\label{subsect;Arf}
Using a certain quadratic form,
Robertello  \cite{Robertello65} introduced 
a knot concordance invariant,
the {\it Arf invariant}  of a knot $K$, 
$\Arf(K)\in \Z_2$. 
Using the Jones polynomial, 
one can define 
this invariant as follows:  
\begin{equation}
V(K;i)=(-1)^{\Arf(K)};
\label{eq;Vi_Arf} \end{equation}
see \cite{Murakami86}. 
By this formula,
we can calculate 
the Arf invariant recursively.

\subsection{The signature} 
\label{subsect;signature}
Murasugi \cite{Murasugi65} introduced 
the \textit{signature} of an oriented link $L$, 
$\Ms(L) \in \Z$.
For a skein triple $(L_+, L_-, L_0)$, Murasugi  showed the follwoing; see  \cite[Lemma~7.1]{Murasugi65}
and \cite[Theorem~1]{Murasugi70}:
\begin{equation} \label{eq;sign_skein2}
\vert \Ms (L_\pm) - \Ms (L_0) \vert \le 1.
\end{equation}
Moreover,
Giller \cite{Giller1982} (cf.~\cite[Theorem 6.4.7]{Murasugi96}) found that
the signature of a knot can be determined 
by the following three axioms:
\begin{itemize}
\item[$\rm(i)$] For the trivial knot $U$, $\Ms(U)=0$.
\item[$\rm(ii)$] Let $(L_+, L_-, L_0)$ be a skein triple. 
If $L_{\pm}$ are knots, 
then 
\begin{equation}
\Ms(L_-)-2 \le \Ms(L_+)\le\Ms(L_-).
\label{eq;sign_skein} \end{equation}
\item[$\rm(iii)$] For a knot $K$, let $\sign V(K; -1)=V(K; -1)/|V(K;-1)|$.  Then
\begin{equation}
(-1)^{\Ms(K)/ 2}=\sign V(K; -1).      
\label{eq;signJones} \end{equation}
\end{itemize}
Note that the signature of a knot $K$, $\sigma(K)$, is  an even integer.


\subsection{The Q polynomial} 
\label{subsect;Qpoly}

\begin{figure}[htbp]
\centerline{\includegraphics[scale=.5]{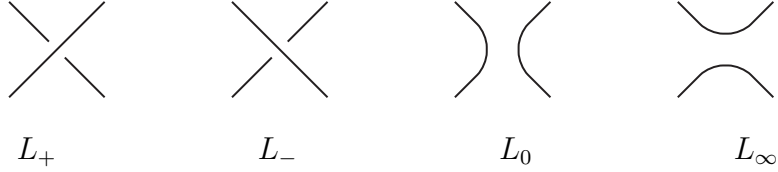} }
\vspace{5pt}\centerline{ $L_{+}$ \hspace{24mm} $L_{-}$ \hspace{24mm} $L_{0}$ 
\hspace{24mm} $L_\infty$\;}
\caption{An unoriented skein quadruple.}
\label{Q_skein}
\end{figure}

The {\it {Q polynomial}}  
$Q(L; z) \in {\Z}[z^{-1}, z]$ \cite{BLM, Ho}
is an invariant of an isotopy type of an unoriented link $L$,
which is defined by the following formulas:
\begin{gather}
 Q(U; z) =1; \\
 Q(L_{+}; z) +Q(L_{-};  z) 
                = z  \left( Q(L_{0};  z) +  Q(L_{\infty};  z) \right).  
                \label{eq;Q_skein}
\end{gather}
where $U$ is the unknot and $L_{+}$, $L_{-}$, $L_{0}$, $L_{\infty}$    
are four unoriented links that are identical except 
near one point where they are as in  Figure~\ref{Q_skein}.
We call $(L_{+}, L_{-}, L_{0}, L_{\infty})$ an {\it unoriented skein quadruple}.

Put $\la(L)=Q\left(L; (\sqrt{5}-1)/2)\right)$.
For a link $L$, Jones \cite{Jones89} has shown:
\begin{equation}
\la(L) = \pm \sqrt{5}^r, 
\label{eq;Jones-Rong_value}\end{equation}
where $r=\dim H_1(\Sigma_2(L); \Z_5)$. 
From  the proof of \cite[Theorem 2]{Rong91}, we have the following:
\begin{proposition}\label{prop;Qpoly}
If two links $L$ and $M$ are related by a single band surgery, then
\begin{equation}
\la(L)/\la(M) \in \{\pm1, \sqrt5^{\pm 1}\}.
\label{eq;prop;Qpoly}\end{equation}
\end{proposition}


This immediately implies the following:
\begin{corollary}\label{cor;Qpoly}  Let $L$ and $M$ be links.
 \par  {\rm(i)}  If $\la(L)/\la(M) =\sqrt5^n$, then $\bd(L,M)\ge |n|$.  
 \par {\rm(ii)} If $\la(L)/\la(M) =-\sqrt5^n$, then $\bd(L,M)\ge |n|+1$.  
\end{corollary}

\begin{remark}\label{remark;Qpoly} 
Since $\bd(J, K)\le \d_2 (J, K)$,
Corollary~\ref{cor;Qpoly} implies 
Corollary~8.1 in \cite{Kanenobu-Miyazawa09} and
Corollary~8.2 in \cite{Kanenobu11}.
\end{remark}


\begin{example} 
(i) Let $K=9_{49}$ or $10_{103}$.  Then $\u_2 (K)=3$ is proved in \cite[p.~453]{Kanenobu-Miyazawa09} by using $\rho(K)=-5$, which further implies $\bu (K)=3$
by Corollary~\ref{cor;Qpoly}(ii).

(ii) Let $F_n$ be the connected sum of $n$ copies of the knot $5_1$.
Since $\la(5_1)=\sqrt 5$ and $\u_2(5_1)=\bu(5_1)=1$, by Corollary~\ref{cor;Qpoly}(i) we have
$\u_2(F_n)=\bu(F_n)=n$.
Since $\la(4_1)=-\sqrt 5$ and $\u_2(4_1)=2$, by Corollary~\ref{cor;Qpoly}(ii) we have
$\u_2(4_1\#F_n)=\bu(4_1\#F_n)=n+2$.
Note that $\u_2 (4_1\# 5_1)=3$ is given in \cite[Table~9.3]{Kanenobu-Miyazawa09}.
\end{example}

\section{Relations between the band- and usual Gordian distances} 
\label{sect;Relations to the band- and usual Gordian distances}
\begin{figure}[hbt]
\centering\includegraphics[width=12cm]{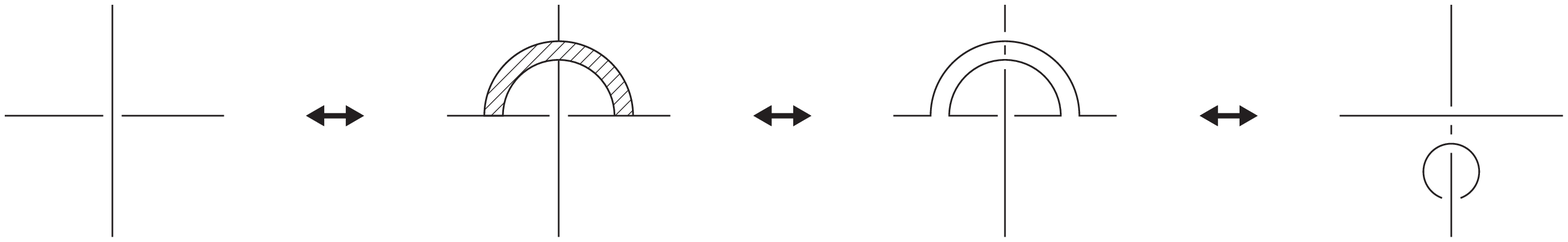}
\vspace{2mm}\par\noindent
\centerline{(a) \hspace{26.1mm} (b) \hspace{26.1mm} (c) \hspace{27.mm} (d)\hspace{.4mm} }
\caption{}\label{F;u1knot}
\end{figure}

We denote by  $\d(J,K)$  the usual Gordian distance between two knots $J$ and $K$.
Then we have the following; see \cite[Theorem~10.1]{Kanenobu11}:
\begin{equation}
\d_2(J,K) \le \d(J,K)+1.
\end{equation}
Then by Theorem~\ref{thm:yasuhara2} we obtain
$\bd(J,K) \le  \d(J,K)+1$.
Moreover, we have the following:
\begin{thm}\label{thm;u1knot} 
For  knots $J$ and $K$,  we have the following: 
\begin{equation}
\bd(J, K)\le
\begin{cases} 
\d(J,K) &{\text{if $\d(J,K)$ is even;}}\\
\d(J,K)+1 &{\text{if $\d(J,K)$ is odd.}}
\end{cases}
\end{equation}
In particular, 
\begin{equation}
\bu(K)\le
\begin{cases} 
\u(K) &{\text{if $\u(K) $ is even;}}\\
\u(K)+1 &{\text{if $\u(K)$ is odd.}}
\end{cases}
\end{equation}
Also, if $J$ and $K$ are unknotting number one knots, then 
\begin{gather}
\bd(J, K)\le 2;
\label{eq;thm;u1knot_1}\\
\bu(J\# K)\le 2.
\label{eq;thm;u1knot_2}\end{gather}
\end{thm}
We prove this using the following Lemma.
\begin{lem}\label{lem;u1knot} 
Suppose that two knots $J$ and $K$ are related by a crossing change,
and let $H$ be the Hopf link.
Then  
\begin{equation}
\bd(J, K\# H)=1.
\label{eq;lem;u1knot_1}
\end{equation}
In particular, if $K$ is an unknotting number one knot.  Then we have:
\begin{gather}
\bd( K, H)=\bu(K\# H)=1.
\label{eq;lem;u1knot_2}
\end{gather}
\end{lem}
\begin{pf}
Suppose that  $J$ is transformed into $K$
by changing the crossing  as shown in Fig.~\ref{F;u1knot}(a).
Attaching a band as shown Fig.~\ref{F;u1knot}(b) to $J$, 
we obtain $K\# H$ as shown in Figs.~\ref{F;u1knot}(c) or (d).
Thus we obtain Eq.~\eqref{eq;lem;u1knot_1}.
%
\end{pf}

\begin{pf*}{Proof of Therem~$\ref{thm;u1knot}$}   If $n=1$, then 
$\bd(J,K)\le \bd(J,K\# H)+\bd(K\# H, K)\le 2$, where we use
Eq.~\eqref{eq;lem;u1knot_2} and the fact that $H$ is transformed into the trivial knot by a single band surgery.
Let $n>1$ and let $J=J_0$, $J_1, \dots, J_n=K$ be a sequence of knots such that $J_{k+1}$ is obtained from $J_k$ by a crossing change.
Then by Lemma~\ref{lem;u1knot}, $J_{k}$ and $J_{k\pm1}\# H$ are related by a band surgery, and so $\bd(J_{k}, J_{k+2})\le 2$.  
Thus we obtain the result.
\end{pf*}
\begin{example} \label{ex;6_1!7_7}
Since $6_1$ and $7_7$ have unknotting number one, 
by  Eq.~\eqref{eq;thm;u1knot_1} in Theorem~\ref{thm;u1knot}
 we have $\bd(6_1!,7_7)\le 2$.
In \cite[Table~5]{Kanenobu11} it is listed that
$\d_2(6_1!,7_7)=2$ or $3$, and so by Theorem~\ref{thm:yasuhara2} we obatin $\bd(6_1!,7_7)= 2$.
\end{example}
\begin{proposition}\label{prop;signature} 
If $J$ and $K$ are knots with $\d_2(J,K)=3$ and $\bu(J,K)=2$,
then $|\Ms(J)-\Ms(K)|\le 2$.
\end{proposition}

\begin{pf} 
By the conditions there exists a 2-component link $L$ such that  $J$ and $L$ are related by an oriented band surgery and that $L$ and $K$ are related by an oriented band surgery. 
Then  by Eq.~\eqref{eq;sign_skein2} we obtain $|\Ms(J)-\Ms(L)|\le 1$ and
$|\Ms(L)-\Ms(K)|\le 1$, which implies the result.
\end{pf}

In \cite[Example~6.4]{Kanenobu11}, several examples of knots with
$H(2)$-unknotting number 3 and pairs of knots with $H(2)$-Gordian distance 3,
which were proved by applying Proposition~\ref{prop;DistanceTwo} below.
In the following, we determine the band-unknotting numbers and band-Gordian distances for these examples.

\begin{figure}[hbt]
\centering\includegraphics[scale=.5]{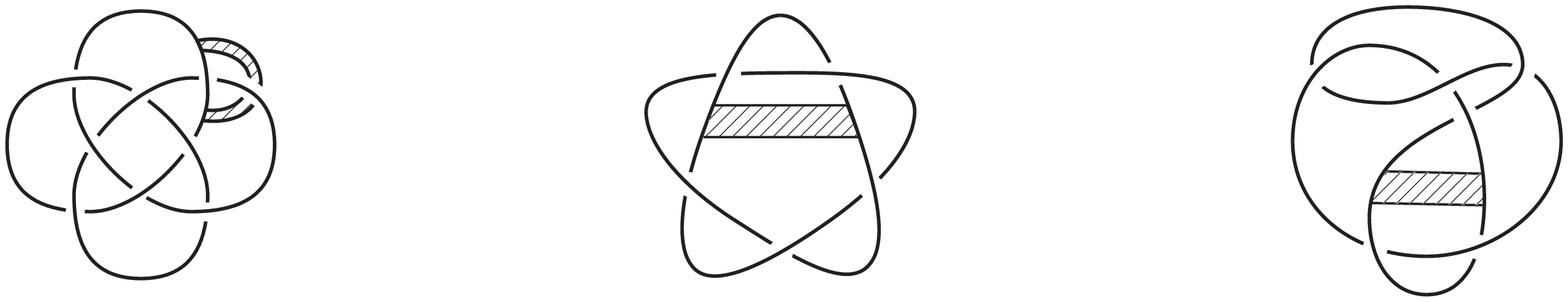}
\centerline{(a) \hspace{42.1mm} (b)  \hspace{42.1mm} (c)}
\caption{}\label{F;u_band(8_18)(5_1)(6_3)}
\end{figure}

\begin{example} \label{ex;8_18} 
(i) In   \cite[Example~6.4(i)]{Kanenobu11} $\u_2 (K)=3$ is proved for $K=8_{18}$, $3_1!\#8_{21}$, $3_1\#9_{40}$, $6_2\#9_{35}$.
Thus by Theorem~\ref{thm:yasuhara2} $\bu(K)=2$  or $3$.
Moreover, we have: 
\begin{equation}
\bu (8_{18})=\bu (3_1!\#8_{21})=2; \quad
\bu (3_1\#9_{40})=\bu (6_2\#9_{35})=3.
\end{equation}
\begin{pf} 
Since $\u(8_{18})=2$, by Theorem~\ref{thm;u1knot} we have  $\bu(8_{18})\le 2$,
and so we obtain $\bu(8_{18})=2$.
Also, attaching the band to the knot $8_{18}$  as shown in Figure~\ref{F;u_band(8_18)(5_1)(6_3)}(a), we obtain the composite link $3_1\# H$, which has band-unknotting number one by Eq.~\eqref{eq;lem;u1knot_2} in  Lemma~\ref{lem;u1knot},  and so $\bu(8_{18})\le 2$.

Since the knots $3_1!$, $8_{21}$ have unknotting number one, by Eq.~\eqref{eq;thm;u1knot_2} in Theorem~\ref{thm;u1knot} we have  $\bu(3_1!\# 8_{21})\le 2$, and so $\bu(3_1!\# 8_{21})=2$.

Since  $\rho(3_1\#9_{40})=-5$, by Corollary~\ref{cor;Qpoly}(ii)
we have  $\bu(3_1\#9_{40})\ge 3$, and so $\bu(3_1\#9_{40})=3$.

Lastly, if $\bu (6_2\#9_{35})=2$, then by Proposition~\ref{prop;signature}
$|\Ms(6_2\#9_{35})|\le 2$. This is a contradiction since $\Ms(6_2)=\Ms(9_{35})=2$, and so
$\Ms(6_2\#9_{35})=4$. Therefore, we obtain  $\bu (6_2\#9_{35})=3$.   
 Notice that since $\Ms(3_1\#9_{40})=4$, we may prove $\bu(3_1\#9_{40})=3$ in the same way.
 \end{pf}

(ii)  In \cite[Example~6.4(ii)]{Kanenobu11},  
$\d_2 (5_1, 3_1\# 3_1)=\d_2 (5_1!, 3_1\# 3_1)=\d_2 (6_3, 3_1\# 3_1!)=3$ are proved.
Thus by Theorem~\ref{thm:yasuhara2} the band-Gordian distance of these pairs are either $2$  or $3$.
Moreover, we have: 
\begin{equation}
\bd (5_1, 3_1\# 3_1)=\bd (6_3, 3_1\# 3_1!)=2; \quad \bd (5_1!, 3_1\# 3_1)=3.
\end{equation}
\begin{pf}
Attaching the band to the knot $5_{1}$ 
as shown in Fig.~\ref{F;u_band(8_18)(5_1)(6_3)}(b), 
we obtain the composite link $3_1\# H$.  
Since $\bd(3_1, H)=1$ 
by Eq.~\eqref{eq;lem;u1knot_2} in Lemma~\ref{lem;u1knot}, 
we have $\bd(3_1\# H, 3_1\# 3_1)=1$, 
and so $\bd(5_1, 3_1\# 3_1) \le 2$.
This also follows from the fact that $\d(5_1, 3_1\# 3_1) = 2$ by using Theorem~\ref{thm;u1knot};  see \cite{DarcySumners}.
Thus we obtain $\bd(5_1, 3_1\# 3_1)=2$.

Attaching the band to the knot $6_{3}$ 
as shown in Fig.~\ref{F;u_band(8_18)(5_1)(6_3)}(c), 
we obtain the composite link $3_1\# H$.  
Since $\bd(3_1!, H)=1$ 
by Eq.~\eqref{eq;lem;u1knot_1} in Lemma~\ref{lem;u1knot}, 
we have $\bd(3_1\# H, 3_1\# 3_1!)=1$, 
and so $\bd(6_3, 3_1\# 3_1!) \le 2$.  Note that $\d(6_3, 3_1\# 3_1!) = 2$; see \cite{DarcySumners}.
Thus we obtain $\bd(6_3, 3_1\# 3_1!)=2$.

We may prove $\bd (5_1!, 3_1\# 3_1)=3$ by Proposition~\ref{prop;signature} since 
$\Ms(5_1!)=-4$ and $\Ms(3_1\# 3_1)=4$.
\end{pf}
\end{example}

\section{Knots with $\bu=2$ and $\u_2=3$}
\label{sect;example}
In this section,  we prove the following theorem. 
\begin{thm}\label{them;main-theorem}
There exist infinitely many knots $K$ such that $\bu(K)=2$ {and} $\u_{2}(K)=3$.  
\end{thm}


We use the following proposition \cite[Theorem 6.3]{Kanenobu11}.
\begin{proposition} 
\label{prop;DistanceTwo}
Let $J$ and $K$ be knots with ${V(J; \om)}/{V(K; \om)} =3$.
If 

{\rm(i)} $\Ms(J)-\Ms(K) \equiv 0  \pmod{8}$, $\Arf(J)\ne\Arf(K)$, or 

{\rm(ii)}
 $\Ms(J)-\Ms(K) \equiv 4  \pmod{8}$,  $\Arf(J)= \Arf(K)$,\\
then
 ${\d_2} (J, K) \ge 3$.
\end{proposition}

\begin{figure}[hbt]
\centering
\includegraphics[scale=.5]{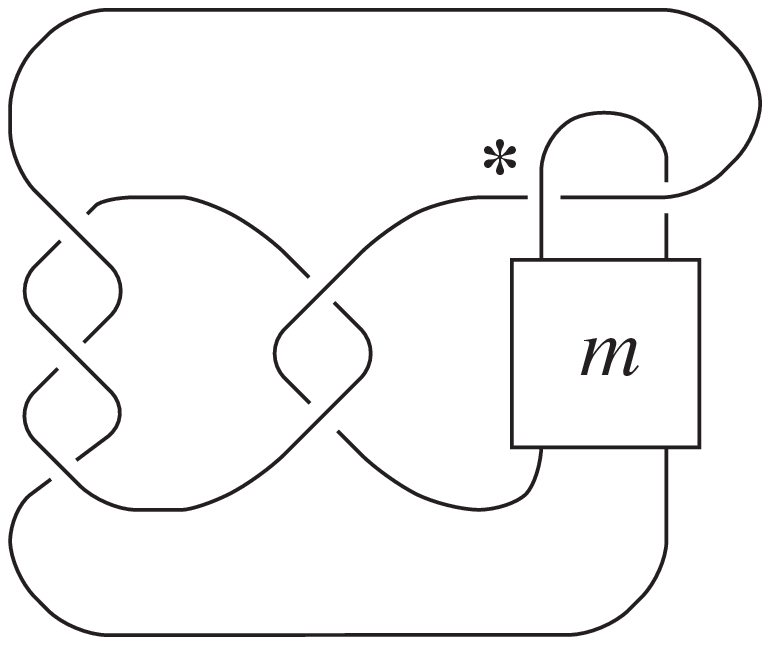}
\caption{}\label{F;knot_Km}
\end{figure}

\begin{figure}[hbt]
\centering
\includegraphics[scale=.5]{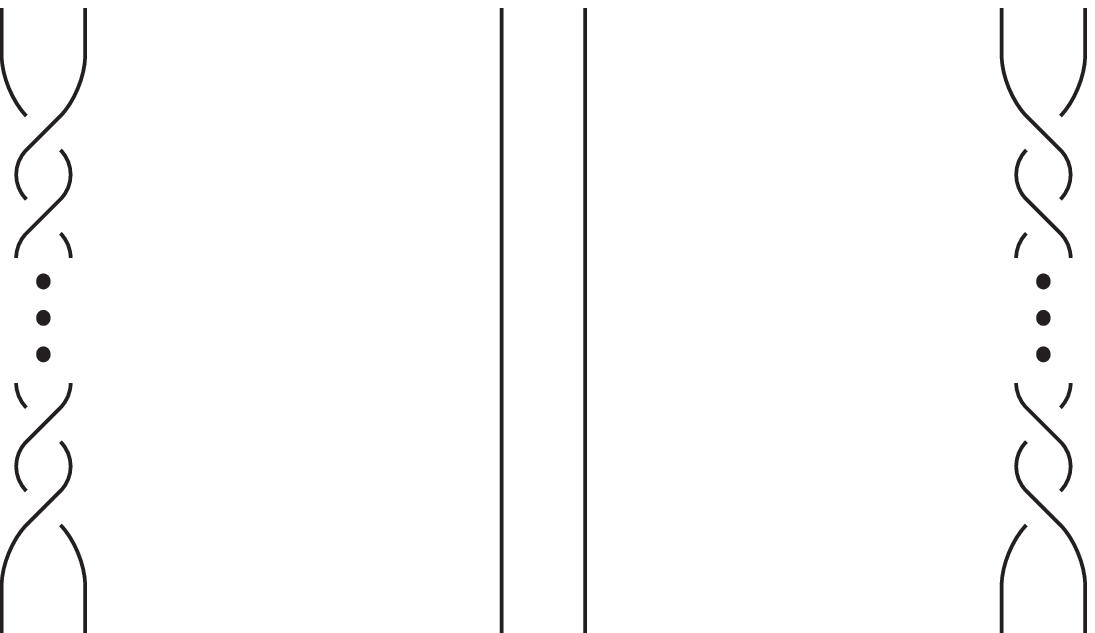}
\vspace{2mm}\par\noindent
\centerline{$m>0$ \hspace{11.5mm} $m=0$ \hspace{11.5mm} $m<0$\hspace{.4mm} }
\caption{}\label{F;2braid}
\end{figure}

Let $K_m$ be a knot as shown in Figure~\ref{F;knot_Km},
where  the tangle labeled $m$ stands for a $2$-braid with $|m|$ crossings in the manner indicated in Fig.~\ref{F;2braid}.
Then we see that 
$K_0=6_3$, 
$K_{-1}=6_2$,
$K_1=8_{21}$, 
$K_{-2}=8_{20}!$,
$K_2=9_{44}$,
$K_{-3}=9_{42}!$,
$K_3=10_{133}$,
$K_{-4}=10_{132}!$. 


We have a skein triple
$(K_{m-2}, K_m, H)$, 
where $H$ is a positive or negative Hopf link 
according as if $m$ is even or odd.
Then from Eq.~\eqref{eq;V_skein}, 
we obtain the following:
\begin{align}
t^{-1} V(K_{m-2};t) - t  V(K_{m};t)  
& = \left( t^{1/2}-t^{-1/2} \right)  V(H;t)  \notag         \\ 
& = 
\begin{cases}
1-t+t^{2}-t^{3}   			&\text{if $m$ is even;}\\
t^{-3}-t^{-2}+t^{-1}-1      &\text{if $m$ is odd;}                        
\end{cases} \label{eq;V_Km}
\end{align}
and 
\begin{align}
V(K_{-1};t)&=V(6_{2};t)=t^{-5}-2t^{-4}+2t^{-3}-2t^{-2}+2t^{-1}-1+t;
\label{eq;V_K-1}\\
V(K_{0};t)&=V(6_{3};t)=-t^{-3}+2t^{-2}-2t^{-1}+3-2t+2t^2-t^3.
\label{eq;V_K0} 
\end{align}

Using these inequalities,
we prove the following three lemmas.
\begin{lem}  \label{lemma:1}
\begin{equation}
\Arf(K_{2n}) = \Arf(K_{2n-1})\equiv   n+1 \pmod2.
\end{equation}
\end{lem}

\begin{pf}
Putting $t=i$ in Eqs.~\eqref{eq;V_Km}, 
\eqref{eq;V_K-1} and \eqref{eq;V_K0}, 
we have $V(K_{n-2}; i)+V(K_{n};i)=0$ and $V(K_{-1};i)=V(K_{0};i)=-1$.
Then using Eq.~\eqref{eq;Vi_Arf}, 
we have
$\Arf(K_{n-2})\ne \Arf(K_n)$ 
and $\Arf(K_0)=\Arf(K_{-1})=1$,
which imply the result.
\end{pf}

\begin{lem}  \label{lemma:2}
\begin{align}
\Ms(K_{2n}) =
\begin{cases}  
 0    &\text{if $n\ge -3$;}\\
-2    &\text{if $n\le -4$;}
\end{cases} 
  \qquad
  \Ms(K_{2n-1}) =
\begin{cases}  
  0    &\text{if $n\le-3$;}\\
    2   &\text{if $n\ge -2$.}
\end{cases} 
\end{align}
\end{lem}

\begin{pf}
The unknotting number of $K_n$ is one. 
Indeed, 
by changing the crossing near the mark $\ast$ 
indicated in Figure~\ref{F;knot_Km}, 
$K_n$ becomes the trivial knot.
Thus by Eq.~\eqref{eq;sign_skein}
and the fact that the signature of a knot is even, 
we have 
\begin{equation}
\Ms(K_n)=0 \text{ or } \pm 2.
\label{eq;signKm1}
\end{equation}
Furthermore, by Eq.~\eqref{eq;sign_skein}, 
we have
\begin{equation}
\Ms(K_{n-2})\le\Ms(K_n).
\label{eq;sign_skeinKm} \end{equation}
Putting $t=-1$ in Eqs.~\eqref{eq;V_Km},
\eqref{eq;V_K-1} and \eqref{eq;V_K0}, 
we have 
\begin{align}
V(K_{n};-1)-V(K_{n-2};-1) & =
\begin{cases}    
-4   \quad &\text{if $n$ is even;}\\
4              &\text{if $n$ is odd;}
\end{cases} 
\end{align}
\[ V(K_{-1};-1)=-11\ \text{and}\ V(K_{0};-1)=13.\]
Thus we obtain the following:
\begin{align}
\sign V(K_{2n};-1)& =
 \begin{cases}    1   \quad &\text{if $n\ge -3$;}\\
                          -1              &\text{if $n\le -4$;}
  \end{cases} 
\\
 \sign V(K_{2n-1};-1) &=
\begin{cases} 
1   \quad &\text{if $n\le-3$;}\\
-1      &\text{if $n\ge -2$.}
\end{cases} 
\end{align}
Therefore,
$\Ms(K_{2n})/2$ is even or odd according as if 
$n\ge -3$ or $n\le -4$; and
$\Ms(K_{2n-1})/2$ is even or odd according 
as if $n\le-3$ or $n\ge -2$.
Then using Eq.~\eqref{eq;sign_skeinKm}, we obtain the result.
\end{pf}

\begin{lem} \label{lemma:3}
\begin{align}
V(K_{2n};\om) &=
 \begin{cases}    
1    &\text{if $n\equiv 0 \pmod3$;}\\
-1           &\text{if $n\equiv 1 \pmod3$;}\\
i\sqrt3      &\text{if $n\equiv 2 \pmod3$;}
\end{cases} 
\\
V(K_{2n-1};\om)& =
\begin{cases}    
1     &\text{if $n\equiv 0 \pmod3$;}\\
- i\sqrt3            &\text{if $n\equiv 1 \pmod3$;}\\
-1            &\text{if $n\equiv 2 \pmod3$.}
\end{cases} 
\end{align}
\end{lem}

\begin{pf}
Putting $t=\om$ in Eqs.~\eqref{eq;V_Km},  \eqref{eq;V_K-1}, and \eqref{eq;V_K0}, 
we have: 
\begin{align}
V(K_{n};\om)+\om V(K_{n-2};\om) & =
\begin{cases}    
\om^2   \quad &\text{if $n$ is even;}\\
-\om^2              &\text{if $n$ is odd;}
\end{cases} 
\end{align}
\[ V(K_{-1};\om)=V(K_{0};\om)=1.\]
Using these inequalities, we obtain the result.
\end{pf}


\begin{pf*}{Proof of Theorem~$\ref{them;main-theorem}$}
Let 
\begin{equation}
J_{l} =
\begin{cases}  
K_{12l+1} \# 3_1!  \quad &\text{if $l\ge 0$;}\\
                          K_{12l-2}  \# 3_1            &\text{if $l\le -1$.}
\end{cases} 
\end{equation}
For each $l \in \Z$ we prove:
\begin{equation}
\bu(J_l)=2,  \quad \u_2 (J_l)=3.      
\label{eq;prop;Jl} \label{prop;Jl}\end{equation}

By lemmas~\ref{lemma:1}, \ref{lemma:2} and \ref{lemma:3}
we obtain Table~\ref{table},
from which we have the following:
\begin{align}
\Ms(J_l)=0, \quad \Arf(J_l)=1, \quad V(J_l;\om)=3,
\end{align}
and so by Proposition~\ref{prop;DistanceTwo}, we have $\u_2 (J_l)\ge 3$.
On the other hand, since  $K_m$ and $3_1$ have unknotting number  one,
by Theorem~\ref{thm;u1knot} we have $\bu(J_l)\le 2$.
Therefore, by Theorem~\ref{thm:yasuhara2}, we obtain $\bu(J_l)= 2$ and $\u_2 (J_l)= 3$.
This completes the proof.
\end{pf*}

\begin{table}[htbp]{\small 
\caption{Unknotting number one knots with $V(\omega)=\pm i\sqrt 3$.} 
\label{table}
\begin{center}
\begin{tabular}{ c|rrrrrrrr}
\noalign{\hrule height0.8pt}
Knots $K$ &$K_{12l+1}$ &$K_{12l+4}$&$K_{12l+7}$&$K_{12l-2}$&$3_1$& $3_1!$ & $6_1$ & $6_1!$ \\
 & $(l\ge 0) $ &  $(l\ge 0) $& $(l\le -1) $& $(l\le -1) $ &&&&\\
\hline
$\sigma(K)$ & 2 &0 & 0& $-2$ &  $2$&$-2$&0&0\\
Arf$(K)$&$0$ &$1$&1&0&1&1&0&0\\
 $V(K;\omega)$&$-i\sqrt 3$&$i\sqrt 3$&$-i\sqrt 3$&$i\sqrt 3$&$-i\sqrt 3$&$i\sqrt 3$&$i\sqrt 3$&$-i\sqrt 3$\\
\noalign{\hrule height0.8pt}
\end{tabular} 
\end{center}
}\end{table}

\begin{remark}
Note that $J_0=8_{21} \# 3_1!$ and $\u_2(J_0)=3$ are
given in Example~6.4(i) in \cite{Kanenobu11}.
\end{remark}

\begin{remark}
Let 
\begin{equation}
\tilde J_{l} =
\begin{cases}  
K_{12l+4} \# 6_1!  \quad &\text{if $l\ge 0$;}\\
                          K_{12l+7}  \# 6_1            &\text{if $l\le -1$.}
\end{cases} 
\end{equation}
Then from Table~\ref{table} we have
$\Ms(\tilde J_l)=0$, $\Arf(\tilde J_l)=1$, $V(\tilde J_l;\om)=3$, and
since $6_1$ has unknotting number one, we obtain
$\bu(\tilde J_l)=2$ and  $\u_2 (\tilde J_l)=3$ for each integer $l$.
\end{remark}

\begin{remark}
Yasuhara  \cite{Yasuhara96} has proved that 
if a knot $K$ bounds a M\"{o}bius band in $B^4$ 
such that $\partial B^4 = S^3$,
then
\begin{equation} \label{eq;Yasuhara}
\Ms (K) - 4\Arf (K) \equiv 0\ \text{ or } \pm 2 \pmod 8.
\end{equation}
If a knot $K$ is transformed into the unknot by a band surgery,
then $K$ bounds a M\"{o}bius band in $B^4$, and so $K$ satisfies
Eq.~\eqref{eq;Yasuhara}, which was generalized in \cite{Kanenobu11},
and Proposition~\ref{prop;DistanceTwo} is a related result.
Another generalization was given by 
Gilmer and Livingston   \cite{Gilmer-Livingston2010}.
\end{remark}

\begin{remark}
Let $C(K)$ be the \textit{crosscap number} of a knot $K$; see \cite{Clark78}. 
Then 
$\d_2 (K) \le C(K)$.
Note that $\bd (K) \le C(K)$ is given in \cite[Prosposition~5.3]{TaniyamaYasuhara94}.
\end{remark}

\section{Band-unknotting numbers and  $H(2)$-unknotting numbers of knots\\ with up to 9 crossings}
\label{section:9cr-knots}

\begin{figure}[htbp]
\centerline{\includegraphics[scale=.5]{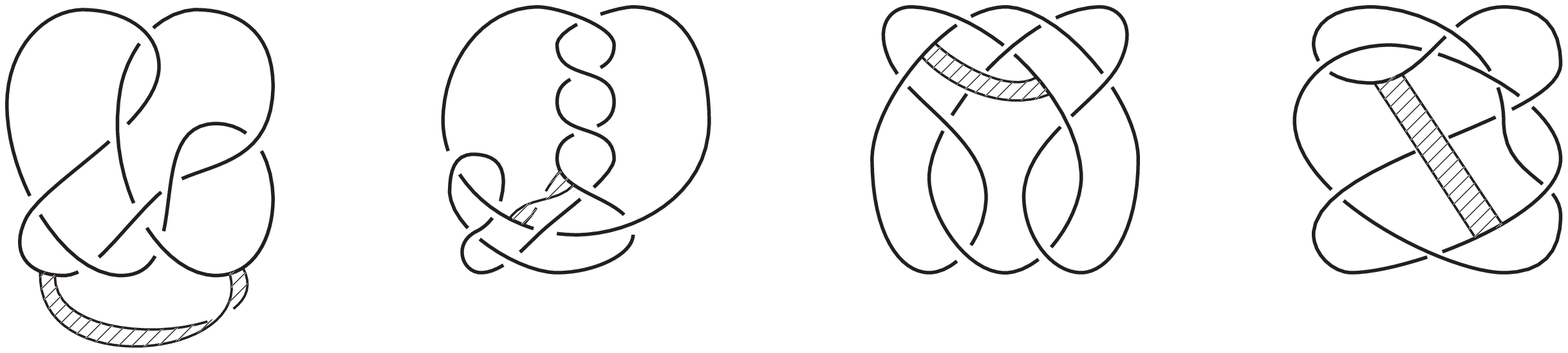}}
\vspace{-5pt}
\centerline{ $8_{14}$ \hspace{27mm} $9_{21}$ \hspace{26mm} $9_{23}$ \hspace{26mm} $9_{26}$ }
\vspace{10pt}
\centerline{\includegraphics[scale=.5]{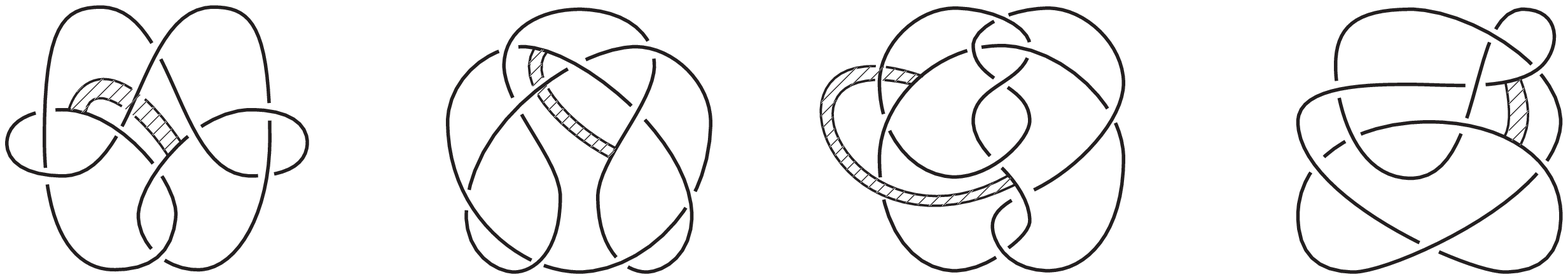} } 
\centerline{$9_{28}$\hspace{26mm} $9_{31}$\hspace{27mm} $9_{32}$ \hspace{26mm} $9_{45}$}
\caption{$H(2)$-unknotting number one knots.}
\label{u2_9cr}
\end{figure}

Table~\ref{table;buH2} lists
the band-unknotting numbers, $H(2)$-unknotting numbers, and usual unknotting numbers
of knots with up to 9 crossings.
In \cite{Kanenobu-Miyazawa09} a table of the $H(2)$-unknotting numbers of knots with up to 9 crossings is given, 
but there were 8 knots whose $H(2)$-unknotting numbers were undecided.
This table is the complete list for the $H(2)$-unknotting numbers of  knots with up to 9 crossings:
In \cite[Example~6.4]{Kanenobu11}  $\u_2(8_{18})=3$ is proved, and
Bao \cite{Bao11} has given a condition for a 2-bridge knot to have $H(2)$-unknotting number one,
which implies the 2-birdge knots $9_{21}$,  $9_{23}$, $9_{26}$, $9_{31}$ have $H(2)$-unknotting number one.
We show these knots and the remaining knots $9_{28}$, $9_{32}$, $9_{45}$ have $H(2)$-unknotting number one in Fig.~\ref{u2_9cr} giving twisted bands that transform into the unknot.

\begin{table}[htbp]{\small 
\caption{Band-unknotting numbers,  $H(2)$-unknotting numbers, and unknotting numbers of knots with up to $9$ crossings.} 
\label{table;buH2}
\begin{center}
\begin{tabular}{ c|c|c|l}
\noalign{\hrule height0.8pt}
 $\bu$ & $\u_2$ & $\u$ & knots \\
\hline
1&1&1&
$3_1$, $6_1$, $6_2$, $7_2$, $7_6$, $8_7$, $8_{11}$, $8_{14}$, $8_{20}$, $9_{19}$,
$9_{21}$, $9_{22}$, $9_{26}$--$9_{28}$, $9_{42}$, $9_{44}$, $9_{45}$.
\\ 
1&1&2&
$5_1$, $7_3$, $7_4$, $8_3$--$8_6$, $8_8$, $8_{10}$, $8_{16}$, $9_4$, $9_5$, $9_7$, $9_8$, $9_{15}$, $9_{17}$, $9_{24}$, $9_{29}$, $9_{31}$, $9_{32}$, $9_{36}$, $9_{43}$.
\\
1&1&3&
$7_1$, $8_{19}$, $9_3$, $9_6$, $9_9$, $9_{13}$.
\\
1&1&4&
$9_1$.
\\
$2$ & $2$ & 1&  $4_1$,   $6_3$, $7_7$, $8_1$, $8_9$,  $8_{13}$, $8_{17}$, $8_{21}$, $9_{2}$,  $9_{12}$, $9_{14}$, $9_{24}$, $9_{30}$, $9_{33}$, $9_{34}$, $9_{39}$.
\\
$2$ & $2$ & 2 &
$7_5$, $8_2$,  $8_{12}$,  $8_{15}$, $9_{11}$, $9_{18}$, $9_{20}$, $9_{37}$,
$9_{40}$, $9_{41}$, $9_{46}$--$9_{48}$, \\
&&&  $3_1\# 3_1$, $3_1!\# 3_1$, $3_1!\# 5_2$, $4_1\# 4_1$,   $3_1\# 6_1$, $3_1!\# 6_1$, $3_1\# 6_2$, $3_1\# 6_3$.
 \\
$2$ & $2$ & 3&
 $9_{10}$, $9_{16}$,  $9_{35}$, $9_{38}$, $3_1\# 5_1$.
\\
$2$ & $2$ & 2-3&
$3_1!\# 5_1$.
\\
$2$ & $3$ & 2&   $8_{18}$.
\\
$3$ & $3$ &  3&  $9_{49}$, $3_1\# 3_1\# 3_1$, $3_1!\# 3_1\# 3_1$, $4_1\# 5_1$.
\\
\noalign{\hrule height0.8pt}
\end{tabular} 
\end{center}
}\end{table}

\begin{remark}
In \cite[Fig.~9.1]{Kanenobu-Miyazawa09} the indicated twisted band  does not transform the knot  $8_{14}$ into the unknot.  The corrected band is shown in Fig.~\ref{u2_9cr}.
In \cite[Table~9.3]{Kanenobu-Miyazawa09} there are errors in the unknotting numbers;
the right ones are:
$\u(3_1\# 5_1)=3$, $\u(3_1!\# 5_1)=2$ or  $3$.
\end{remark}

\section{Special value of the Jones polynomial and a slice knot} 
\label{section:slice-knots}
\label{sect:Special values of the Jones polynomial and a slice knot} 

A knot in $S^3$ is \textit{slice} if it bounds a disk in the 4-ball $B^4$ such that
$\partial B^4 =S^3$.
Using an oriented version of Theorem~\ref{thm;Jpoly} we have proved  \cite[Corollary~4.5]{Kanenobu10}:
If $K$ is a ribbon knot, then  $V(K; \om) \neq -1$.
 Moreover, we have:
\begin{proposition} \label{prop:slice}
If $K$ is a  slice knot, then $V(K; \om) \neq -1$.
\end{proposition}

\begin{pf}
%
Let $K$ be a slice knot with $V(K; \om) = -1$.
First, note that by Eq.~\eqref{eq;LickorishMillett86} the determinant of $K$ is not  divided by 3.
There are polynomial $P(t)$ in $\Z[t^{-1}, t]$ and integers $a$, $b$ such that 
\begin{equation}
V(K;t)=(t^2-t+1)P(t)+at+b.
\label{eq;pf_prop:slice1}\end{equation}
Then 
we have  $a=0$ and $b=-1$, and so we obtain:
\begin{equation}
V(K;-1)=3P(-1)-1.
\label{eq;pf_prop:slice2}\end{equation}
On the other hand, the value $|V(K;-1)|$ is the determinant of $K$, which is a square integer.  More precisely, $V(K;-1)=\Delta_K(-1)$,
the value of the Conway-normalized Alexander polynomial  $\Delta_K(t)$ at $t=-1$, where $\Delta_K(t)$ satisfies $\Delta_K(t)=\Delta_K(t^{-1})$ and $\Delta_K(1)=1$;  see \cite[Chapter~16]{Lickorish97}.
Since $K$ is a slice knot,  $\Delta_K(t)$ has the form $\Delta_K(t)=f(t)f(t^{-1})$, where $f(t)$ is a polynomial in $\Z[t]$ \cite{MR0211392}; 
cf.~\cite[Theorem~8E20]{Rolfsen}. 
Putting $f(-1)=3k\pm1$, $V(K;-1)=9k^2\pm 6k +1$, which contradicts 
Eq.~\eqref{eq;pf_prop:slice2}.  This completes the proof.
\end{pf}

\section*{Acknowledgments}
The authors would like to thank Akira Yasuhara who suggested Corollay~\ref{cor:yasuhara2}.
The first author was partially supported by KAKENHI,
Grant-in-Aid for Research Activity start-up (No. 00614009),
Japan Society for the Promotion of Science.
The second author was partially supported by KAKENHI, Grant-in-Aid for Scientific Research (C) (No. 21540092), Japan Society for the Promotion of Science.


\end{document}